\documentclass[a4paper,11pt]{amsart}

\usepackage[applemac]{inputenc}
\usepackage{amsmath, amsthm, amssymb}
\usepackage{eucal, enumerate, cite, graphicx}
\usepackage[all]{xy}

\theoremstyle{theorem}
\newtheorem{theorem}{Theorem}[section]
\newtheorem{proposition}[theorem]{Proposition}
\newtheorem{corollary}[theorem]{Corollary}
\newtheorem{lemma}[theorem]{Lemma}
\newtheorem{conjecture}[theorem]{Conjecture}
\theoremstyle{definition}

\theoremstyle{remark}
\newtheorem{remark}[theorem]{Remark}

\newcommand{\N}{\mathbb{N}}
\newcommand{\Z}{\mathbb{Z}}
\newcommand{\Q}{\mathbb{Q}}
\newcommand{\R}{\mathbb{R}}

\renewcommand{\P}{\mathbb{P}}

\DeclareMathOperator{\Nef}{Nef}
\DeclareMathOperator{\Psef}{Psef}

\begin{document}
\bibliographystyle{alpha}

\title[Rational curves on Calabi-Yau threefolds]{On a conjecture of Oguiso about rational curves on Calabi-Yau threefolds}

\author{Simone Diverio \and Andrea Ferretti}
\address{Simone Diverio \\ Institut de Math\'ematiques de Jussieu - CNRS\\
Universit\'e Pierre et Marie Curie \\ Paris}
\email{diverio@math.jussieu.fr} 
\address{Andrea Ferretti \\ Laboratoire Paul Painlevé \\ Université Lille 1 \\ Lille}
\email{ferrettiandrea@gmail.com}

\dedicatory{Dedicated to the memory of Marco Brunella}

\date{\today}

\begin{abstract}
Let $X$ be a Calabi-Yau threefold. We show that if there exists on $X$ a non-zero nef non-ample divisor then $X$ contains a rational curve, provided its second Betti number is greater than $4$. 
\end{abstract}
\maketitle
\section{Introduction}

Let $X$ be a Calabi-Yau threefold, that is, a compact complex Kähler manifold of dimension three with trivial canonical bundle $K_X \simeq \mathcal O_X$ and finite fundamental group. 
 
In a series of foundational papers starting with \cite{Wil89}, Wilson began a systematic study of the geometry of Calabi-Yau threefolds by looking at the structure of their ample cone and deducing from that several remarkable kind of algebraic fiber space structures on them (especially in the case of large Picard number). In particular, he obtained as a consequence that every Calabi-Yau threefold whose Picard number is greater than $19$ always contains a rational curve (this result was later improved to Picard number greater than $13$ in \cite{HB-W92}). 
 
Using the same circle of ideas, Peternell showed in \cite{Pet91} (see also \cite{Ogu93}) that, under the condition of the existence of a non-zero effective non-ample divisor, one can recover a rational curve on $X$; roughly speaking, this is because such a divisor defines a fibration which is not an isomorphism so that one tries to get some positive dimensional fiber in which one hopes to find the desired curve. Of course, the existence of such a divisor forces the Picard number of $X$ to be greater than one; moreover, one can suppose that every such effective divisor is nef, otherwise by the Cone Theorem one would get immediately a rational curve. Thus, one should ask what happens if one merely has a non-zero nef non-ample divisor; this is the content of a conjecture proposed by Oguiso:

\begin{conjecture}[Oguiso \cite{Ogu93}]\label{conj}
Let $X$ be a Calabi-Yau threefold possessing a non-zero nef non-ample divisor $D$. Then, $X$ contains a rational curve.
\end{conjecture}

In \cite{Ogu93} and independently in \cite{Wil94}, some partial answers are given, assuming for instance that the numerical dimension of $D$ is one, or that $D$ intersects non-trivially the second Chern class of $X$. Here we give a positive answer to Conjecture \ref{conj} with just an extra (mild) hypothesis on the Picard number (or, equivalently on the second Betti number) of $X$.

\begin{theorem}
Let $X$ be a Calabi-Yau threefold. If there exists on $X$ a non-zero nef non-ample divisor, then $X$ does contain a rational curve, provided its second Betti number is greater than $4$. 
\end{theorem}

We also have some partial result -- depending on the existence of special divisors on $X$ satisfying certain numerical conditions -- for the remaining cases $b_2(X) = 2, 3, 4$. Such results can be summarized in the following:

\begin{proposition}
Let $X$ be a Calabi-Yau threefold possessing a non-zero nef non-ample divisor $D$. Then, $X$ contains a rational curve provided one of the following conditions is fulfilled:
\begin{itemize}
\item[(i)] the null cone $\mathcal N_X$ of $X$ is irreducible and $b_2(X)=4$,
\item[(ii)] the null cone $\mathcal N_X$ of $X$ is irreducible, $b_2(X)=3$ and there do not exist two $\Q$-divisors $E,F$ such that
$$
D^2\cdot E=0,\quad E^3=0,\quad E^2\cdot F=0,\quad E\cdot F^2=0,
$$
and $D,E,F$ span $NS(X)_{\R}$.
\item[(iii)] $b_2(X)=2$ and either there exists a divisor $E$, not a multiple of $D$, such that $E^3=0$, or there exists a $\Q$-divisor $E$, not a multiple of $D$, and such that $3(E^2\cdot D)^2=4E^3\times(E\cdot D^2)$.
\end{itemize}
\end{proposition}

The null cone of $X$ is defined as the locus in $NS(X)_{\R}$ of divisors whose top self-intersection vanishes: being a polynomial locus in $NS(X)_{\R}$, its irreducibility in the statement is intended as an algebraic variety. Moreover, we shall see that on a Calabi-Yau threefold without rational curves all nef non-ample divisors must sit in it (cf. Proposition \ref{nefpsef}).

\bigskip

One of our motivations to look at this kind of problem was to get some evidence in the direction of a conjecture by Kobayashi, which states that every projective Kobayashi hyperbolic manifold has ample canonical bundle. It turns out (see Section 4 for more details) that in order to prove this conjecture in dimension three it suffices to show that Calabi-Yau threefolds are not Kobayashi hyperbolic, \textsl{i.e.} they do admit a non-constant entire (\textsl{a priori} transcendental) map from the complex plane. Of course, such a map exists if the manifold contains a rational or an elliptic curve (or more generally a non-constant holomorphic image of a complex torus). Another perhaps more appropriate way to show the non-hyperbolicity of $X$ would be to exhibit a sequence of curves of general type $(C_\ell)$ in $X$ such that
$$
-\frac{\chi(\widehat C_\ell)}{\deg C_\ell}\to 0\quad\text{as}\quad \ell\to\infty,
$$
where $\widehat C_\ell$ is the normalization of $C_\ell$ and the degree is taken with respect to any polarization of $X$ -- see \cite{Dem97} for more details.

The results stated above thus permit, as in \cite{Pet91}, to exclude a certain number of cases to be checked in order to prove such a non-hyperbolicity statement. Finally, as far as we know, there is no known example of a Calabi-Yau threefold without either rational curves, elliptic curves or non-constant holomorphic images of complex tori.
 
\subsection{Notation and conventions}

By a Calabi-Yau threefold we shall always mean a compact Kähler manifold of dimension three with trivial canonical bundle and finite fundamental group. In particular, this implies that the irregularity $q(X)=\dim H^1(X,\mathcal O_X)$ is zero. Since on a compact complex threefold
\begin{equation*}
H^1(X,\mathcal O_X) \cong \bigl(H^2(X,K_X)\bigr)^*
\end{equation*}
by Serre's duality, a Calabi-Yau threefold $X$ satisfies $H^2(X,\mathcal O_X)=0$ and thus it is always projective. From $H^1(X,\mathcal O_X)=H^2(X,\mathcal O_X)=0$, one easily deduces that $\operatorname{Pic}(X)$ is isomorphic to $H^2(X,\Z)$ and that they are both isomorphic to the Néron-Severi group of $X$. In particular, the Picard number of $X$ equals its second Betti number. 

With our definition, one also has $c_2(X)\ne 0$. To see this, just recall (see for instance \cite{Kob87}) that a compact Kähler manifold such that $c_1(X)=c_2(X)=0$ in $H^\bullet(X,\mathbb R)$ is a finite unramified quotient of a torus. Thus, the fundamental group of $X$ must contain a free abelian group of rank $6$, which is impossible by our assumption on $\pi_1(X)$.

For the definitions, basic properties and notations about positive cones and positivity concepts for divisors we refer to \cite{Laz04}, while for the birational geometry of pairs and related classical results we follow \cite{K-M98}.

Finally, by a rational point, we shall always mean a $\Q$-rational point.

\subsubsection*{Acknowledgments} We are glad to thank F.~Campana, B.~Claudon and J.-P.~Demailly for their valuable comments on a preliminary version of this paper. Also, we would like to thank P. M. H. Wilson for pointing out to us the paper \cite{Wil94} and E. Brugallé for stimulating discussions about real cubic forms.

\section{Structure of the nef cone and rational curves}

Let $X$ be a projective manifold. As usual, define the nef cone $\Nef(X)$ and the pseudoeffective cone $\Psef(X)$ of $X$ to be the closed convex cones in the real Néron-Severi space generated respectively by the classes of nef and effective divisors on $X$; it always holds $\Nef(X)\subseteq\Psef(X)$ and the interior of these two cones gives the open convex cones respectively of ample and of big classes. On $NS(X)_{\R}$ it is defined a (integral) top-intersection form
$$
NS(X)_{\R}\ni D\mapsto D^{\dim X}\in\R.
$$
Its zero locus will be denoted by $\mathcal N_X\subseteq NS(X)_{\R}$ and is usually called the \emph{null cone}. Finally, the \emph{nef boundary} $\mathcal B_X\subseteq NS(X)_{\R}$ is the boundary $\partial\Nef(X)$ of the nef cone.

\begin{proposition}\label{nefpsef}
Let $X$ be a projective manifold with trivial canonical bundle and no rational curves. Then, $\Nef(X)=\Psef(X)$ and moreover the nef boundary $\mathcal B_X$ is entirely contained in the null cone $\mathcal N_X$.
\end{proposition}

\begin{proof}
Let $D$ be a big $\Q$-divisor -- in particular $D$ is $\Q$-effective. By taking a small multiple of $D$, we can suppose that the pair $(X,D)$ is Kawamata log terminal (klt). Since $K_X\simeq\mathcal O_X$, then $K_X+D=D$ and the Cone Theorem tells us immediately that $D$ is nef, otherwise we would have some negative extremal ray generated by the class of a rational curve in $X$. Thus, the interior of $\Psef(X)$ is contained in $\Nef(X)$ and therefore $\Psef(X)\subseteq\Nef(X)$. 

The second assertion is quite general and holds for any (smooth) projective manifold such that $\Nef(X)=\Psef(X)$. Since the nef cone is contained in (one connected component of) the locus $\{D^{\dim X}\ge 0\}$, if $\mathcal B_X$ was not contained in $\mathcal N_X$ we would find a nef $\R$-divisor $D$ with strictly positive top self-intersection $D^n>0$. Then, $D$ is a big $\mathbb R$-divisor, that is, its class lies in the interior of $\Psef(X)$. But in this case it cannot be on $\mathcal B_X$.
\end{proof}

Thus, on a projective manifold with trivial canonical bundle and no rational curves, every effective divisor is nef. Now, let $D$ be a non-zero effective $\Q$-divisor on $X$. Similarly as above, by taking a small (rational) multiple of $D$, we can suppose that the pair $(X,D)$ is log canonical (lc). Since $D=K_X+D$ is nef, by the log-abundance for threefolds \cite{K-M-MK94,K-M-MK04}, the linear system associated to some multiple of $D$ is free, so that it defines an algebraic fiber space structure $\phi_D\colon X\to Y$, with $Y$ a normal projective variety whose dimension
\begin{equation*}
\dim Y = \kappa(X,D)
\end{equation*}
equals the Kodaira dimension of $D$ (which coincides, since $D$ is abundant, with its numerical dimension $\nu(X,D)=\max\{k\in\N\mid D^k\not\equiv 0\}$) and $D\sim_\Q\phi^*_DA$ for some ample line bundle $A$ on $Y$. 

\begin{theorem}[Peternell \cite{Pet91}, see also \cite{Ogu93}]
Let $X$ be a smooth projective threefold with trivial canonical bundle and $c_2(X)\ne 0$, and let $D$ be as above. Assume that $D$ is not ample, so that $D\in\mathcal B_X$; then $X$ does contain a rational curve.
\end{theorem}

Remark that the hypothesis $c_2(X)\ne 0$ is necessary: otherwise $X$ is a finite unramified quotient of a torus and therefore it cannot contain any rational curve. We shall sketch and slightly rephrase the proof here below, since it clarifies the general strategy and fixes notations for the entire business.

\begin{proof}[Outline of the proof (see \cite{Ogu93} for more details)]
To begin with, observe that since $c_2(X)\ne 0$, by the Beauville-Bogomolov decomposition theorem (see also Section 4), after a finite étale cover, $X$ is either a Calabi-Yau threefold or a product of a projective $K3$ surface and of an elliptic curve. In the latter case, by \cite{M-M83}, we find a rational curve on the $K3$ factor and hence on $X$. So, from now on, we can suppose that $X$ is a Calabi-Yau threefold.

The idea is to find a rational curve in the fibers of the fibration $\phi_D\colon X\to Y$ defined by $D$. By Proposition \ref{nefpsef}, we can suppose that every big divisor is ample. Thus $\kappa(X,D)=\nu(X,D)$ is either $1$ or $2$ and the fibers of $\phi_D$ have codimension $\nu(X,D)$ in $X$. Moreover, by adjunction, the (very) general fiber of $\phi_D$ has trivial canonical bundle (since $D\sim_\Q\phi^*_DA$ is trivial when restricted to fibers). 

If $\nu(X,D)=2$, then it can be shown that $Y=W$ is a rational surface and the general fiber of $\phi_D$ is a smooth elliptic curve. A finer analysis shows that, depending on $c_2(X)\cdot D$ the following two cases can occur:
\begin{itemize}
\item[($II_+$)] $c_2(X)\cdot D>0$; then $W$ has only quotient singularities and $\phi_D$ is an elliptic fibration with at least one singular fiber and with no multiple fibers; this gives us the desired rational curve as a singular fiber.
\item[($II_0$)] $c_2(X)\cdot D=0$; then $W$ is non-Gorenstein with only quotient singularities and $\phi_D$ is a smooth elliptic fibration in codimension one over $W$; moreover there is a non-Gorenstein point $w\in W$ such that $\dim\phi_D^{-1}(w)=2$ so that $\phi_D$ is never equi-dimensional.
\end{itemize}
In order to find a rational curve in case of fibration of type ($II_0$), one can then use the following result, which we state in a simplified form that suffices for our aims:

\begin{theorem}[Kawamata \cite{Kaw91}]\label{exc}
Let $f\colon X\to Y$ be a surjective projective morphism, where $X$ is smooth and $-K_X$ is $f$-nef (that is, it intersects non negatively the curves which are contracted by $f$). Then, any irreducible component of $\operatorname{Exc}(f)=\{x\in X\mid \dim  f^{-1}(f(x))>\dim X-\dim Y\}$ is uniruled.
\end{theorem}

Suppose now $\nu(X,D)=1$. In this case, $Y=\mathbb P^1$ and the general fiber of $\phi_D$ is either a $K3$ surface or an abelian surface. The first case corresponds to $c_2(X)\cdot D>0$ and by the existence of rational curves on any projective $K3$ surface \cite{M-M83}, we are done. If $c_2(X)\cdot D=0$, then nothing can be said using merely the semiample fibration associated to $D$.

In this case, starting from $D$, we try to construct another divisor $N$ for which one of the two following statements applies.
\begin{lemma}[Key Lemma, \cite{Wil89}]\label{keylemma}
Let $X$ be a Calabi-Yau threefold. Assume that there exists an ample divisor $H$ and a non-nef divisor $N$ on $X$ such that
$$
N^3>0,\quad N^2\cdot H>0\quad\text{and}\quad N\cdot H^2>0.
$$
Then, $X$ contains a rational curve.
\end{lemma}

\begin{proposition}[Oguiso \cite{Ogu93}]\label{semiample}
Let $X$ be a Calabi-Yau threefold and $N$ a non-zero nef divisor on $X$ such that $c_2(X)\cdot N>0$. Then $N$ is $\Q$-effective (and hence semiample).
\end{proposition}

In order to do that, consider the affine line (of rational slope) of divisors $N_t=H-tD$, where $H$ is an ample divisor and $t\in\Q$, and the inequalities (suggested by Lemma \ref{keylemma}):
$$
N_t^3>0,\quad N_t^2\cdot H>0\quad\text{and}\quad N_t\cdot H^2>0.
$$
Since $c_2(X)\ne 0$ and $X$ is supposed to be not uniruled, then by \cite{Miy87} we have $c_2(X)\cdot\Nef(X)\ge 0$ and the inequality is strict on ample classes. From $D^2\equiv 0$, one can deduce that the previous system of inequalities is equivalent to
\begin{equation}\label{keyinequalities}
t<t_0=\frac{H^3}{3D\cdot H^2},
\end{equation}
and moreover $t_0$ is such that $N_{t_0}^3=0$ and $N_{t_0}^2\cdot H>0$. 

If $N_t$ is not nef for some solution $t$ of (\ref{keyinequalities}), then we are done by applying Lemma \ref{keylemma} to $N=N_t$ (perhaps after perturbing a little bit $t$ in order to deal with rational divisors). Otherwise, $N_t$ is nef for every rational solution of (\ref{keyinequalities}), therefore $N_{t_0}$ is also nef (and $t_0$ is of course rational). Since $c_2(X)\cdot N_{t_0}=c_2(X)\cdot H>0$, Proposition \ref{semiample} with $N=N_{t_0}$ gives effectiveness of $N_{t_0}$. Finally, $N_{t_0}^2\cdot H>0$ implies $\nu(X,N_{t_0})=2$ and therefore $N_{t_0}$ endows $X$ with the structure of an algebraic fiber space of type ($II_+$).
\end{proof}

Observe that in the last part of the proof, in order to construct the new divisor $N$, we did not use any effectivity property of $D$, but just its nefness and the fact that its numerical dimension was one. As a byproduct, one obtains:

\begin{proposition}[Oguiso \cite{Ogu93}, see also \cite{Wil94}]\label{c2}
Let $X$ be a Calabi-Yau threefold possessing a non-zero nef $\Q$-divisor $D\in\mathcal B_X$ such that either $c_2(X)\cdot D\ne 0$ or $\nu(X,D)=1$. Then, $X$ contains a rational curve.
\end{proposition}

Let us rephrase the statement in the following form, more convenient for what follows.

\begin{proposition}\label{summa}
Let $X$ be a smooth projective threefold with trivial canonical bundle, $c_2(X)\ne 0$ and without rational curves. Then, $X$ is a Calabi-Yau threefold and every rational divisor $D\in\mathcal B_X$ is such that $\nu(X,D)=2$ and $c_2(X)\cdot D=0$.
\end{proposition}

\section{Arithmetic of the nef boundary}

Using the description of the nef boundary given in the previous section, we now prove the following result.

\begin{theorem} \label{main}
Let $X$ be a smooth projective Calabi-Yau threefold with $b_2(X) > 4$, and assume there exists a non-zero rational divisor $D$ on the nef boundary of $X$. Then, $X$ does contain a rational curve.
\end{theorem}

The \textsl{naïve} idea of the proof is as follows. Start with the given divisor $D$ which is by hypothesis rational, non-zero, and on the nef boundary. By Proposition \ref{summa}, we can assume that $\nu(X,D)=2$ and $c_2(X)\cdot D=0$, otherwise we are done. By \cite{Miy87}, $c_2(X)\cdot A>0$ for every ample divisor $A$ on $X$, so that the hyperplane in $NS(X)_{\R}$ defined by $c_2(X)\cdot x=0$ is a supporting hyperplane for the convex cone $\Nef(X)$ at $D$. By Proposition \ref{nefpsef}, we may assume that $\mathcal B_X$ is locally described by the cubic equation $x^3=0$, so that the tangent space at the point $D$ is given by $\{D^2\cdot x=0\}$. On the other hand, $\nu(X,D)=2$ tells us that this is a true hyperplane, that is, $\mathcal B_X$ is smooth at $D$, and thus the supporting hyperplane $\{c_2(X)\cdot x=0\}$ is in fact the tangent hyperplane to $\mathcal B_X$ at $D$ and coincides with the hyperplane $\{D^2\cdot x=0\}$.

Now, two different scenarios are possible: either $\mathcal N_X$ is irreducible or it factorizes as a union of the hyperplane $\{c_2(X)\cdot x=0\}$ and a residual quadratic locus (we shall see that it can never factorize as a product of three linear loci as soon as $b_2(X)\ge 4$). By Proposition \ref{c2}, it will suffice to find a non-zero rational point $E\in\mathcal B_X$ which satisfies $c_2(X)\cdot E\ne 0$. In the first case, this will be achieved by finding such an $E$ by an iteration argument and, in the second case, by finding it as a non-zero rational point sitting on the residual quadric (and not on the hyperplane).

\bigskip

\begin{figure}[htbp]
\centering
\includegraphics{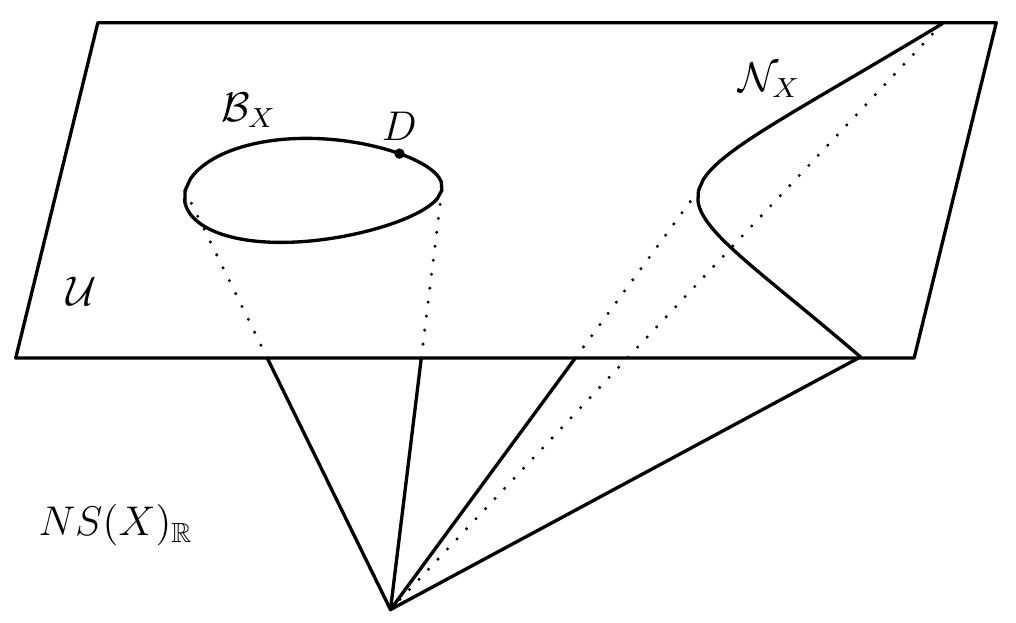}
\caption{The null cone and the nef boundary cut by the affine hyperplane $\mathcal U$ at the rational nef non-ample divisor $D$.}\label{fig:cones}
\end{figure}

Let us now make some general remarks and state some elementary lemmata before entering into the proof. First of all, it is quite natural to look at things in the projective space $\P(NS(X)_{\R})$. We shall fix once and for all an affine open chart $\mathcal U$ of $\P(NS(X)_{\R})$ containing all directions of $\Nef(X)$: one has to think at it as an affine hyperplane of rational slope in $NS(X)_{\R}$ passing through $D$ and cutting $\Nef(X)$ transversally (see \figurename~\ref{fig:cones} for a qualitative picture of the situation); with a slight abuse of language we shall continue to call $\mathcal N_X$ and $\mathcal B_X$ their projectivized in $\P(NS(X)_{\R})$ and talk about convexity and boundedness properties referring to the open affine part $\mathcal U$. 

From now on, we will work in this real projective space of real dimension $b_2(X)-1$. Notice that $\mathcal B_X$ bounds a compact convex set in $\mathcal U$.

\begin{remark}\label{segment}
Let $\Omega\subset\R^n$, $n\ge 2$, an open bounded set, $x\in\partial\Omega$ smooth point of the boundary and $T_{\partial\Omega,x}$ the tangent space of $\partial\Omega$ at $x$. If $\ell$ is any line passing through $x$ and not contained in $T_{\partial\Omega,x}$, then $\ell$ must meet $\partial \Omega$ in at least one more point different from $x$ and, of course, not contained in $T_{\partial\Omega,x}$.
\end{remark}

\begin{lemma}\label{singquadric}
Let $Q \subset \R^n$, $n\ge 2$, be a singular quadric and $H \subset \R^n$ an affine hyperplane such that $\complement(Q\cup H)$ has a bounded connected component $K$. Then the singular set of $Q$ consists of a single point which lies on $\partial K\setminus H$.
\end{lemma}

In particular, three affine hyperplanes can never bound a compact set in $\R^n$, as soon as $n\ge 3$.

\begin{proof}

To prove that $Q$ has just one singular point it is enough to prove that $Q \cap H$ is a smooth quadric. Indeed the singular locus of $Q$ is an affine space, hence it is a point if it is disjoint from $H$. Notice that $Q \cap H$ bounds a compact set in $H$; since a singular quadric is a cone, it cannot bound a compact set, hence $Q \cap H$ is smooth.

It remains to check that the vertex of $Q$ lies on $\partial K$. If this is not the case, let $x \in \partial K$ be any point outside $H$; then the half-line starting from $x$ and going away from the vertex must lie on $\partial K$, contradicting its boundedness.
\end{proof}

\begin{lemma}\label{triple}
Let $P \subset \R^n$, $n\ge 3$, be an irreducible cubic and $x \in P$ be a smooth point with tangent space $T = T_{P,x}$. Suppose that $P|_{T}$ is everywhere non-reduced. Then, for every neighbourhood $U$ of $x$, $P\cap U$ cannot lie on only one side of $T$. In particular, $P$ cannot bound a convex set near $x$.
\end{lemma}

\begin{proof}
Observe that a cubic which is everywhere non-reduced is necessarily given by a linear equation raised to the cube. We can suppose without loss of generality that $x = 0 \in \R^n$ and, after a linear transformation, that $T = \{x_n=0\}$. Moreover, we can arrange the coordinates so that $P$ is given by the equation
$$
p(x) = x_1^3 + x_n\,q(x),
$$
where $q$ is a quadratic polynomial. For every $\varepsilon \in \R$ consider the equation
\begin{equation*}
\phi_{\varepsilon}(t) = t^3 + \varepsilon q(t, 0, \dots, 0, \varepsilon) = 0.
\end{equation*}
Since it is cubic, there is always at least one solution. For $\varepsilon = 0$ we have the triple solution $t = 0$; hence every zero of $\phi_{\varepsilon}$ remains small for $\varepsilon$ small. In particular there are arbitrarily small points on $P$ both in the case $\varepsilon > 0$ and $\varepsilon < 0$, which is the thesis.
\end{proof}

\begin{proof}[Proof of Theorem \ref{main}]
Suppose first that the cubic locus $\mathcal N_X$ is reducible. Since $b_2(X)\ge 5$, by Lemma \ref{singquadric} it cannot be the union of three hyperplanes. Therefore $\mathcal N_X$ splits as a union $\mathcal H\cup\mathcal Q$, where $\mathcal H$ is a hyperplane and $\mathcal Q$ the residual irreducible (indefinite) quadric. Since $\mathcal N_X$ is defined over $\Q$, then $\mathcal H$ and $\mathcal Q$ are defined over $\Q$, too; this can be seen straightforwardly by acting on the defining equations with automorphisms of $\mathbb C$ which by hypothesis leave $\mathcal N_X$ invariant. Observe moreover that in this case $\Nef(X)$ is bounded by $\mathcal H\cup\mathcal Q$.

Of course, rational points are dense in $\mathcal H$; by Meyer's theorem \cite{Ser77}, since $b_2(X)\ge 5$, $\mathcal Q$ satisfies the  Hasse-Minkowski local-global principle and thus it does have a rational point. If this point is smooth, then rational points of $\mathcal Q$ are dense by projection and thus rational points are also dense in $\mathcal N_X$ (and hence in $\mathcal B_X$). So, we find (plenty of) rational divisors in $\mathcal B_X$ whose intersection with $c_2(X)$ is non-zero, and we are done thanks to Proposition \ref{c2}. If the rational point is singular for $\mathcal Q$, then by Lemma \ref{singquadric} it lies in $\mathcal B_X$. Thus, we have found a rational point on the nef boundary with numerical dimension $1$ (otherwise it would be a regular point) and we conclude again by Proposition \ref{c2}. 

Next, suppose that $\mathcal N_X$ is irreducible. In this case the hypothesis on the second Betti number can be slightly weakened to $b_2(X)\ge 4$. Let $D$ be the (smooth) rational point of $\mathcal B_X$ as in the hypotheses and consider the tangent space $T_{\mathcal B_X,D}=T_{\mathcal N_X,D}$ of $\mathcal N_X$ at $D$: as we have seen, it is given by $\{D^2\cdot x=0\}=\{c_2(X)\cdot x=0\}$, and hence it is defined over $\Q$. The intersection 
$$
C=\mathcal N_X\cap T_{\mathcal N_X,D}
$$
is a cubic of dimension one less, which is by construction singular at $D$, hence rational over $\Q$.  Unfortunately, this is not enough to find rational divisors on the boundary of the ample cone which are not on $T_{\mathcal N_X,D}=\{c_2(X)\cdot x=0\}$. Denote by
$$
S := C(\Q)
$$
the set of rational points thus obtained.
If one of the points in $S$ is singular for $\mathcal N_X$, we are done as above by projection; otherwise we can repeat the procedure starting from any point $D' \in S$ and produce more rational points on $\mathcal N_X$.

We claim that we can suppose that repeating the procedure will actually yield some point not in $S$: thus we obtain a rational point $E$ of $\mathcal N_X$ such that either $E$ is already contained in $\mathcal B_X$ and $c_2(X)\cdot E\ne 0$ or $E\in\mathcal N_X\setminus\mathcal B_X$ and the line through $D$ and $E$ is not contained in $T_{\mathcal N_X,D}$. In the first case we are done by Proposition \ref{c2}. In the second case the conclusion follows from Remark \ref{segment} and the fact that if a line meets a cubic defined over $\Q$ at three points, two of which are rational, then the third one is rational as well.

To prove the claim, suppose that the tangent space to $\mathcal N_X$ at all points in $S$ is the same. Since the cubic $C$ is rational over $\Q$, rational points on $C$ are dense in the Zariski topology, and it follows that the tangent space to $\mathcal N_X$ at all points of $C$ is the same. This implies that $C$ is everywhere singular, hence everywhere non-reduced. But then, by Lemma \ref{triple}, $\mathcal B_X$ cannot be convex at the point $D$, contradiction.
\end{proof}

As a byproduct of the proof above, we obtain immediately the following corollaries, which improve slightly two results contained in \cite{Wil89}.

\begin{corollary}
Let $X$ be a Calabi-Yau threefold whose second Betti number $b_2(X)>4$. If $\mathcal N_X$ is reducible, then $X$ contains a rational curve.
\end{corollary}

Remark that in the above corollary, there is no assumption on the existence of a \lq\lq special\rq\rq{} divisor on $X$; 	unluckily, we are unable to give any satisfactory description of when such a situation actually occurs.

\begin{corollary}
Let $X$ be a Calabi-Yau threefold whose second Betti number $b_2(X)>3$. If $\mathcal N_X$ is irreducible and there exists a non-zero nef non-ample $\Q$-divisor on $X$, then $X$ contains a rational curve.
\end{corollary}

Let us now state a partial answer to Oguiso's conjecture in the case when $b_2(X)=3$.

\begin{proposition}\label{b2=3}
Let $X$ be a Calabi-Yau threefold with a non-zero nef non-ample $\Q$-divisor $D$. Suppose that $\mathcal N_X$ is irreducible, $b_2(X)=3$ and there do not exist two $\Q$-divisors $E,F$ such that
\begin{equation}\label{inflection}
D^2\cdot E=0,\quad E^3=0,\quad E^2\cdot F=0,\quad E\cdot F^2=0,
\end{equation}
and $D,E,F$ span $NS(X)_{\R}$. Then, $X$ contains a rational curve.
\end{proposition}

Unfortunately, if $\mathcal N_X$ is reducible, it seems to us  that nothing can be really said as soon as $b_2(X)\le 4$.

\begin{proof}
By hypothesis, $\mathcal N_X$ is a cubic curve in $\P(NS(X)_{\R})\simeq\P^2(\R)$, which is smooth at $D$. By choosing a plane which cuts $\Nef(X)$ transversely, we will look at $\mathcal N_X$ as an affine cubic and talk about bounded and unbounded components, as before (see again \figurename~\ref{fig:cones}).

First, we can assume that $\mathcal N_X$ is smooth. Otherwise its singular locus is defined over $\Q$; since it consists of only one single point, that point is rational. But then the cubic is rational over $\Q$ and we are done.

By Harnack's theorem, $\mathcal N_X$ must have at most two real connected components, and since part of $\mathcal N_X$ bounds the nef cone, it must have exactly two components, one of which bounded. 
Take the tangent line to $\mathcal N_X$ at $D$: its equation is $\{D^2\cdot x=0\}$ (which can be assumed to coincide with $\{c_2(X)\cdot x=0\}$) and it meets $\mathcal N_X$ at one more rational (smooth) point $E$ which satisfies then $E^3=D^2\cdot E=0$ and lies on the unbounded component; clearly, the tangent line to $\mathcal N_X$ at $E$ cannot be the same as the one tangent at $D$. 

Now, suppose that $E$ is not an inflection point for $\mathcal N_X$. Then, repeating the same construction we find another rational point $D'$ on (the unbounded component of) $\mathcal N_X$ such that the line through it and $D$ is different from the tangent line to $\mathcal N_X$ at $D$. We conclude then by Remark \ref{segment} since the third point obtained must be rational.

In remains to show that $E$ is not an inflection point of $\mathcal N_X$, but a straightforward computation shows that this is the case exactly when one can find another divisor $F$ satisfying the last two conditions in (\ref{inflection}) and such that $D,E,F$ span $NS(X)_{\R}$. 
\end{proof}

To finish with, let us make a final remark concerning Oguiso's conjecture in the case $b_2(X)=2$, which is the smallest possible dimension of $H^2(X,\R)$ in order to have a non-zero nef non-ample divisor. In such a situation, the nef cone is bounded by two extremal rays, one of which is rational by hypothesis. Call $D$ a generator of this rational ray: we can suppose that it satisfies $c_2(X)\cdot D=0$ and $D^2\not\equiv 0$, otherwise we find as usual a rational curve in $X$. Thus, in order to find a rational curve by these methods one has to show that the other extremal ray is rational, too. Fix another rational divisor $E$ linearly independent with $D$: the null cone $\mathcal N_X$ with respect to this basis is described by
\begin{equation}\label{nullb2=2}
E^3\,x^3+3E^2\cdot D\,x^2y+3E\cdot D^2\,xy^2=0.
\end{equation}
Therefore, its (projective) solutions other than $D=[0:1]$ are rational if and only if
$$
9(E^2\cdot D)^2-12E^3\times(E\cdot D^2)
$$
is a perfect square. This is the case if one can choose either $E$ to be in $\mathcal N_X$ (but not necessarily in the nef boundary $\mathcal B_X$) or such that $3(E^2\cdot D)^2=4E^3\times(E\cdot D^2)$ (and $E^3\ne 0$). The first possibility tells that we can thus find a non-zero nef non ample divisor intersecting non-trivially $c_2(X)$. The latter possibility means that (\ref{nullb2=2}) has, apart from $D$, only one other solution of multiplicity two given by $[-3E^2\cdot D:2E^3]$, that is, the other extremal ray is spanned by the rational divisor 
$$
D':=(-3E^2\cdot D)E+(2E^3)D,
$$
whose numerical dimension is $1$ (by direct computation of $(D')^2\cdot E=(D')^2\cdot D=0$). By Proposition \ref{c2}, this discussion gives:

\begin{proposition}\label{b2=2}
Let $X$ be a Calabi-Yau threefold with a non-zero nef non-ample $\Q$-divisor $D$ and suppose that $b_2(X)=2$. Then $X$ contains a rational curve provided 
\begin{itemize}
\item either there exists a $\Q$-divisor $E$, not a multiple of $D$, such that $E^3=0$,
\item or there exists a $\Q$-divisor $E$, not a multiple of $D$, and such that $3(E^2\cdot D)^2=4E^3\times(E\cdot D^2)$.
\end{itemize}
\end{proposition}

\section{On a conjecture of Kobayashi}

Let $X$ be a compact complex space. Recall that $X$ is Kobayashi hyperbolic if and only if there is no non-constant entire holomorphic map $f\colon\mathbb C\to X$. In particular, if $X$ is Kobayashi hyperbolic, then it cannot contain any non-hyperbolic subvariety.

In 1970, S. Kobayashi proposed the following:

\begin{conjecture}\label{kob}
Let $X$ be a smooth complex projective manifold. If $X$ is Kobayashi hyperbolic then the canonical bundle $K_X$ is ample.
\end{conjecture}

We take the opportunity here to reproduce a quite standard argument in order to reduce this conjecture to showing that projective varieties $X$ with non-positive Kodaira dimension $\kappa(X)\le 0$ are not hyperbolic. Observe that if $X$ is uniruled (and hence certainly not Kobayashi hyperbolic) then $\kappa(X)=-\infty$; on the other hand, it is a conjecture that $\kappa(X)=-\infty$ implies uniruledness. This conjecture is known to hold true in dimension less than or equal to three.

So, let $X$ be a smooth Kobayashi hyperbolic projective manifold. First of all, by the celebrated criterion of Mori, $K_X$ is nef -- otherwise $X$ would contain a rational curve. Thus, $K_X$ is already in the closure of the ample cone (it is somewhat surprising that Kobayashi could formulate Conjecture \ref{kob} without Mori's criterion at his disposal). 

Suppose for a moment that one can prove that Kobayashi hyperbolicity implies strictly positive Kodaira dimension. Then, we can apply the following to $K_X$.

\begin{theorem}[Iitaka fibrations]
Let $X$ be a normal projective variety and $L\to X$ a line bundle on $X$ such that $\kappa(X,L)>0$. Then, for all sufficiently large $k$ such that $H^0(X,kL)\ne 0$, the rational mappings $\phi_k\colon X\dashrightarrow Y_k$ induced by the linear system $|kL|$ are birationally equivalent to a fixed algebraic fiber space
$$
\phi_\infty\colon X_\infty\to Y_\infty
$$
of normal varieties, and the restriction of $L$ to a very general fiber of $\phi_\infty$ has zero Kodaira-Iitaka dimension.

More specifically, there exists for any such large $k$ a commutative diagram 
$$
\xymatrix{
X \ar@{-->}[d]_{\phi_k} & X_\infty \ar[l]_{u_\infty} \ar[d]^{\phi_\infty}\\
Y_k & Y_\infty \ar@{-->}[l]^{\nu_k}
}
$$
of rational maps and morphisms, where the horizontal maps are birational and $u_\infty$ is a morphism. One has that $\dim Y_\infty=\kappa(X,L)$ and moreover, if we set $L_\infty=u_\infty^*L$ and take $F\subset X_\infty$ to be a very general fiber of $\phi_\infty$, then
$$
\kappa(F,L_\infty|_F)=0.
$$
\end{theorem}
The Iitaka fibration of an irreducible variety $X$ is by definition the Iitaka fibration associated to the canonical bundle (on any non-singular model) of $X$, provided $\kappa(X)>0$. A very general fiber $F$ of the Iitaka fibration of $X$ satisfies $\kappa(F)=0$. 

We claim that this implies that the Kodaira dimension of a projective Kobayashi hyperbolic manifold $X$ must be maximal, that is, $X$ is of general type. Indeed, if $1\le\kappa(X)<\dim X$, then $\phi_\infty$ has positive dimensional fibers and the very general ones have zero Kodaira dimension. Since $u_\infty$ is birational,  the (Zariski closure of the) image by $u_\infty$ of a very general fiber of $\phi_\infty$ will have zero Kodaira dimension, so that we would, by our assumptions, find a non-Kobayashi hyperbolic positive dimensional subvariety of $X$, contradiction.

Now, if $K_X$ is big and nef, the Base Point Free Theorem implies immediately that it is semi-ample, so that a large multiple of $K_X$ defines a genuine surjective, generically $1$-$1$ morphism $f\colon X\to X'$ whose exceptional locus $\operatorname{Exc}(f)$ is empty if and only if $K_X$ is ample; by Theorem \ref{exc}, this locus, if non-empty, has uniruled irreducible components. Thus, if $X$ is hyperbolic and $K_X$ is big (and nef, by Mori), then $K_X$ must be ample (see also \cite{Tak08}).

Next, how to prove that a projective manifold $X$ of non-positive Kodaira dimension is not hyperbolic? We restrict our attention to the three-dimensional case (for some partial result on hyperKähler manifold in higher dimension, one can see also \cite{Cam92}), since in dimension two the conjecture follows from the birational classification plus the existence of rational curves on projective $K3$ surfaces \cite{M-M83}. As observed above, we already know that a projective threefold $X$ of negative Kodaira dimension is uniruled, so that we can suppose $\kappa(X)=0$, that is, $h^0(X,mK_X)\in\{0,1\}$ for all $m\ge 0$. But then, since the abundance conjecture holds true in dimension three, a large multiple of $K_X$ must be globally generated by only one global section; in other words some large multiple of $K_X$ is trivial.

Thus, $K_X$ is torsion and $c_1(X)\in H^2(X,\mathbb R)$ is zero. By the Beauville-Bogomolov decomposition theorem \cite{Bea83}, a manifold with vanishing real first Chern class is, up to finite étale covers, a product of complex tori, Calabi-Yau manifolds and holomorphic symplectic manifolds. Hence at least one factor of $X$ is a torus, or $X$ is a Calabi-Yau. Since complex tori are obviously not Kobayashi hyperbolic, one is reduced to showing that Calabi-Yau threefolds are not Kobayashi hyperbolic (since Kobayashi hyperbolicity is preserved under étale covers).

It is now clear that results that suggest that Calabi-Yau threefolds are not hyperbolic also go in the direction of proving the Kobayashi conjecture. In particular, if the Kobayashi conjecture is false, there must exist a hyperbolic Calabi-Yau threefold; by our results together with \cite{Pet91,HB-W92,Ogu93}, such a variety $X$ will have the following properties:
\begin{enumerate}[(i)]
\item
the second Betti number $b_2(X) \leq 13$;
\item
every non-zero effective divisor on $X$ is ample;
\item
if $b_2(X) \geq 5$, every non-zero nef divisor on $X$ is ample. 
\end{enumerate}
Of course, for $b_2(X) \leq 4$, we have given a partial description in Propositions \ref{b2=3} and \ref{b2=2}.

Needless to say, our results deal with the existence of rational curves while to prove non-hyperbolicity of Calabi-Yau threefolds one could look at the much weaker condition of the existence of entire (\textsl{a priori} just transcendental) curves: hopefully this will be the subject of further investigations, possibly via the study of Ahlfors' currents and their cohomology classes.

\bibliography{bibliografia}{}

\end{document}